\magnification\magstep1
\baselineskip = 18pt
\def\n{\noindent}
\null\vskip.2in
\input mssymb.tex

\centerline{{\bf The plank problem for symmetric bodies}\footnote{}{AMS
1980 Subject classification: \ 52A37, 46A22\hfil\break $^{(1)}$Partially
supported by NSF DMS-8807243}}

\centerline{by}
\centerline{Keith Ball$^{(1)}$}\bigskip\medskip

\centerline{Department of Mathematics}
\centerline{Texas A\&M University}
\centerline{College Station, TX \ 77843}\vskip.4in

\n {\bf Abstract.} Given a symmetric convex body $C$ and $n$ hyperplanes in
an Euclidean space, there is a translate of a multiple of $C$, at least
${1\over n+1}$ times as large, inside $C$, whose interior does not meet any
of the hyperplanes. The result generalizes Bang's solution of the plank
problem of Tarski and has applications to Diophantine approximation.
\vfill\eject

\n {\bf \S 1. Introduction and preliminary observations.}

In the 1930's, Tarski posed what came to be known as the plank problem. A
{\it plank} in ${\Bbb R}^d$ is the region between two distinct parallel
hyperplanes. Tarski conjectured that if a convex body of minimum width $w$
is covered by a collection of planks in ${\Bbb R}^d$, then the sum of the
widths of these planks is at least $w$. Tarski himself proved this for the
disc in ${\Bbb R}^2$. The problem was solved in general by Bang [B]. At the
end of his paper, Bang asked whether his theroem could be strengthened by
asking that the width of each plank should be measured relative to the
width of the convex body being covered, in the direction of the normal to
the plank. This affine invariant plank problem has a number of natural
formulations: \ in particular, as the multi-dimensional ``pigeon-hole
principle'' stated in the abstract. The history of the affine plank problem
from Bang's paper to the present, together with many interesting remarks
can be found in the papers [Gr], [R] and especially [Ga].

In the case of symmetric bodies, the problem is perhaps most naturally
stated in terms of normed spaces. Let $X$ be a normed space. A plank in
$X$, is a region of the form

$$\{x\in X\colon \ |\phi(x) - m|\le w\}$$

\n where $\phi$ is a functional in $X^*, m$ a real number and $w$ a
positive number. If $\phi$ is taken to be a functional of norm 1, $w$ is
said to be the {\it half-width} of the plank. The theorem proved here is
the following.	  \medskip

\n {\bf Theorem 1.} If the unit ball of a Banach space $X$ is covered by a
(countable) collection of planks in $X$, then the sum of the half-widths of
these planks is at least 1.

The theorem is obviously best possible in the sense that for every unit
vector $\phi \in X^*$, the ball of $X$ can be covered by one (or more)
planks, perpendicular to $\phi$, whose half-widths add up to 1.

The infinite-dimensional case of Theorem~1 does not follow formally from
the finite-dimensional: \ it will be discussed and proved in Section~3 of
the paper. For finite-dimensional spaces, one can restate Theorem~1, with
the aid of  compactness, as follows.\medskip

\n {\bf Theorem 2.} If $(\phi_i)^n_1$ is a sequence of unit functionals on
a finite-dimensional normed space $X, (m_i)^n_1$ is a sequence of reals
and $(w_i)^n_1$ a sequence of positive numbers with $\sum w_i = 1$ then
there is a point $x$ in the unit ball of $X$ for which

$$|\phi_i(x) - m_i|\ge w_i \quad \hbox{for every}\quad i.$$

The question answered by Theorem~2 arises quite naturally in the theory of
badly approximable numbers. In his paper [D], Davenport made use of the
following observation. If $C$ is a cube in ${\Bbb R}^d$ and $(H_i)^n_1$ are
$n$ hyperplanes, then there is a cube $C'$ at least $2^{-n}$ times as
large
as $C$, inside $C$, with faces parallel to those of $C$, whose interior is
not met by any $H_i$. This pigeon-hole principle can be strengthened
considerably if Theorem~2 is invoked. (This was already noticed by
Alexander in [A].) The estimate below immediately transfers to give
sharper estimates in Davenport's theorems.\medskip

\n {\bf Corollary.} If $C$ is a convex body, with a center of symmetry, in
${\Bbb R}^d$ and $(H_i)^n_1$ are hyperplanes, then there is a set of the
form $x + {1\over n+1} C$ inside $C$, whose interior is not met by any
$H_i$. The result is obviously sharp for every $n$ and $C$.

\n {\bf Proof.} Assume that $C$ is centered at the origin and let $X$ be
the normed space represented on ${\Bbb R}^d$ with unit ball $C$. For each
$i$, choose a functional $\phi_i$ of norm 1 in $X^*$ and a real number
$m_i$ so that

$$H_i = \{x\in {\Bbb R}^d\colon \ \phi_i(x) = m_i\}.$$

\n By Theorem~2 there is a point $x\in {n\over n+1} C$ with

$$|\phi_i(x) - m_i|\ge {1\over n+1} \quad {\rm for \ each}\quad  i.$$

\n Then the set

$$x + {1\over n+1} C \subset C$$

\n and for every $y$ in $x + {1\over n+1} C, \|y-x\| \le {1\over n+1}$ so
that for each $i, |\phi_i(y) - \phi_i(x)| \le {1\over n+1}\colon$ \ hence
$\phi_i(y) - m_i$ has the same sign as $\phi_i(x) - m_i$. Thus, for each
$i$, the whole of $x + {1\over n+1}C$ lies on the same side of $H_i$ as
$x$ does.$\hfill \square$\medskip

Theorem~2 is readily reduced to a combinatorial theorem concerning
matrices. For a sequence $(\phi_i)^n_1$ of norm 1 functionals on $X$,
construct a matrix $A = (a_{ij})$ given by

$$a_{ij} = \phi_i(x_j),\qquad 1\le i,j\le n$$

\n where for each $j, x_j$ is a point in the unit ball of $X$ at which
$\phi_j$ attains its norm; i.e. $\phi_j(x_j) = \|x_j\| = 1$. If
$(\lambda_j)^n_1$ is a sequence of reals with

$$\sum|\lambda_j|\le 1$$

\n then the vector $x = \sum \lambda_j x_j$ has norm at most 1 and for each
$i$,

$$\phi_i(x) = \sum_j a_{ij}\lambda_j.$$

\n Thus, Theorem~2 follows from\medskip

\n {\bf Theorem 2$'$.} Let $A = (a_{ij})$ be an $n\times n$ matrix whose
diagonal entries equal 1, $(m_i)^n_1$ a sequence of reals and $(w_i)^n_1$ a
sequence of non-negative numbers with $\sum\limits ^n_i w_i\le 1$. Then
there is a sequence $(\lambda_j)^n_1$ with

$$\sum_j |\lambda_j|\le 1$$

\n and, for each $i$,

$$\left|\sum_j a_{ij}\lambda_j-m_i\right| \ge w_i.$$

It is also easy to see that Theorem~2$'$ follows immediately from
Theorem~2 by regarding the rows of such a matrix as unit vectors in
$\ell^n_\infty$. Bang effectively proved Theorem~2$'$ for symmetric
matrices: \ his elegant argument is reproduced here as a lemma, since the
precise statement will be needed later.\medskip

\n {\bf Lemma 3}  (Bang). Let $H = (h_{ij})$ be a real, symmetric $n\times
n$ matrix with 1's on the diagonal, $(\mu_i)^n_1$ a sequence of reals and
$(\theta_i)^n_1$ a sequence of non-negative numbers. Then there is a
sequence of signs $(\varepsilon_j)^n_1 \ (\varepsilon_j = \pm 1$ for each
$j$) so that for each $i$,

$$\left|\sum_j h_{ij}\varepsilon_j\theta_j-\mu_i\right| \ge \theta_i.$$

\n {\bf Proof.} Choose signs $(\varepsilon_j)^n_1$ so as to maximise

$$\sum_{i,j} h_{ij} \varepsilon_i\varepsilon_j \theta_i\theta_j - 2 \sum_i
\varepsilon_i \theta_i\mu_i.$$

\n Fix $k \ (1\le k\le n)$ and define $(\delta_j)^n_1$ by

$$\delta_j = \left\{\matrix{\hfill\varepsilon_j&{\rm if}&j\ne k\hfill\cr
&&\cr
\hfill-\varepsilon_j&{\rm if}&j=k.\hfill\cr}\right.$$

\n Then

$$0\le \sum_{i,j} h_{ij}\varepsilon_i\varepsilon_j\theta_i\theta_j - 2
\sum_i \varepsilon _i\theta_i\mu_i - \left( \sum_{i,j} h_{ij}
\delta_i\delta_j \theta_i\theta_j - 2 \sum_i \delta_i\theta_i
\mu_i\right)$$

\n and since $H$ is symmetric this expression is

$$4\varepsilon_k\theta_k \sum_{j\ne k} h_{kj}\varepsilon_j\theta_j-4
\varepsilon_k\theta_k\mu_k.$$

\n Since $h_{kk} = 1$, the latter is

$$-4 \theta^2_k + 4\varepsilon_k \theta_k \sum_j h_{kj}\varepsilon_j
\theta_j - 4\varepsilon_k \theta_k\mu_k,$$

\n and so

$$\eqalign{4\theta^2_k &\le 4\varepsilon_k\theta_k \left(\sum_j h_{kj}
\varepsilon_j\theta_j-\mu_k\right)\cr
&\le 4 \theta_k\left| \sum_j h_{kj}\varepsilon_j \theta_j -
\mu_k\right|.}$$

\n Since this holds for each $k$, the result is proved.$\hfill \square$
\medskip

Note that the hypothesis of symmetry cannot be dropped from Lemma~3: \
consider, for example, the matrix $\left(\matrix{\hfill 1&1\cr \hfill
-1&1}\right)$ for $\theta_1 = \theta_2 = 1$ and $\mu_1=\mu_2 = 0$.

In the proof of Theorem~2$'$ it may be assumed that $w_i = {1\over n}, \ 1
\le i\le n$ since planks of varying widths can be almost covered by
slightly overlapping ``sheets'', all of the same width. (This ``change of
density'' argument is not really needed but simplifies the succeeding
arguments.) It will be shown that in this situation, Theorem~2$'$ can be
strengthened: \ the sequence $(\lambda_j)^n_1$ to be chosen will actually
satisfy

$$\sum_j \lambda^2_j \le {1\over n}$$

\n (which implies that $\sum\limits_j |\lambda_j|\le 1$ by the
Cauchy-Schwartz inequality). This stronger statement can be attacked by
Hilbert space methods: \ if the satement is true for $AU$ where $U$ is an
orthogonal matrix, then it is true for $A$.  Unfortunately, it is not the
case that for every matrix $A$ with 1's on the diagonal, there is an
orthogonal matrix $U$ with $AU$ both symmetric and having large diagonal.
For example, if

$$A = \left(\matrix{1&1\cr {1\over 2}&1\cr}\right)$$

\n then the only symmetric matrices of the form $AU$ are

$$\pm\left(\matrix{1&1\cr 1&{1\over 2}\cr}\right) \quad {\rm and}\quad
\pm\left(\matrix{{5\over \sqrt{17}}\hfill&{3\over \sqrt{17}}\hfill\cr
{3\over \sqrt{17}}\hfill&{7\over 2\sqrt{17}}\hfill\cr}\right).$$

\n Nevertheless, Theorem~2$'$ is proved by using a modified matrix $A$.
\bigskip

\n {\bf \S 2. Symmetrisations of matrices and the proof of the main
theorem.}

The modification of a matrix, needed for the proof, is described by the
following lemma. From now on, if $H$ is a matrix, $H$ will be said to be
positive if it is symmetric and positive semi-definite.
\medskip

\n {\bf Lemma 4.} Let $A$ be an $n\times n$ matrix of reals, each of whose
rows is non-null. Then there is a sequence $(\theta_i)^n_1$ of positive
numbers and an orthogonal matrix $U$ so that the matrix $H = (h_{ij})$
given by

$$h_{ij} = \theta_i(AU)_{ij}$$

\n is positive and has 1's on the diagonal.

Lemma~4 can be proved using a fixed point theorem or other topological
methods. However it has an elementary proof which provides an alternative
description of the sequence $(\theta_i)^n_1$. Recall that for a matrix $B$,
the trace-class, or nuclear, norm $\|B\|_{C_1}$ of $B$, is $tr(H)$, where
$H$ is the positive square root of $BB^*$. By the Cauchy-Schwartz
inequality

$$\|B\|_{C_1} = \max\{tr(BU)\colon \ U \ {\rm orthogonal}\}.$$

\n Also by the Cauchy-Schwartz inequality, if $B$ and $C$ are $n\times n$
matrices then

$$\|BC\|_{C_1} \le (tr(BB^*))^{1/2} (tr(CC^*))^{1/2}.$$

Before the proof of Lemma~4 it will be convenient to prove the lemma that
really forms the crux of the proof of Theorem~2. The estimate is somewhat
unusual since it involves the sum of the squares of the diagonal entries of
a matrix: \ nevertheless, it is a consequence of the matrix Cauchy-Schwartz
inequality.\medskip

\n {\bf Lemma 5.} If $H = (h_{ij})$ is a positive matrix with non-zero
diagonal entries and $U$ is orthogonal then

$$\sum_i {(HU)^2_{ii}\over h_{ii}} \le \sum_i h_{ii}.$$

\n {\bf Proof.} For each $i$ let $\gamma_i = {(HU)_{ii}\over h_{ii}}$ and
let $D$ be the diagonal matrix, $diag(\gamma_i)^n_1$ and $T$, the positive
square root of $H$. Then

$$\eqalign{\sum_i {(HU)^2_{ii}\over h_{ii}} &= \sum_i \gamma_i(HU)_{ii}\cr
&= tr \ DHU \le \|DH\|_{C_1}\cr
&= \|(DT)T\|_{C_1}\cr
&\le [tr \ DT(DT)^*]^{1\over 2} [tr \ TT^*]^{1\over 2}\cr
&= [tr \ DHD]^{1\over 2} [tr \ H]^{1\over 2}\cr
&= \left(\sum_i \gamma^2_i h_{ii}\right)^{1\over 2} \left(\sum_i
h_{ii}\right)^{1\over 2}\cr
&= \left(\sum {(HU)^2_{ii}\over h_{ii}}\right)^{1\over 2} \left(\sum
h_{ii}\right)^{1\over 2}.}$$

\n Hence

$$\left(\sum_i {(HU)^2_{ii}\over h_{ii}}\right)^{1\over 2} \le \left(\sum_i
h_{ii}\right)^{1\over 2}.\eqno \square$$

\n  Lemma 5 immediately implies:

\n {\bf Lemma 6.} If $H = (h_{ij})$ is a positive $n\times n$ matrix with
non-zero diagonal entries, then

$$\left\|\left( {1\over \sqrt{h_{ii}}} h_{ij}\right)\right\|_{C_1} \le
\sqrt{n} \|H\|^{1/2}_{C_1}.$$

\n {\bf Proof.} There is some orthogonal matrix $U$ for which

$$\left\|\left( {1\over \sqrt{h_{ii}}} h_{ij}\right)\right\| _{C_1} =
\sum_i {1\over \sqrt{h_{ii}}} (HU)_{ii}.$$

\n By the Cauchy-Schwartz inequality this is at most

$$\sqrt{n} \left(\sum_i {(HU)^2_{ii}\over h_{ii}}\right)^{1\over 2}$$

\n and this is at most $\sqrt{n} \big(\sum h_{ii}\big)^{1\over 2} =
\sqrt{n} \|H\|_{C_1}$ by Lemma~5.$\hfill \square$
\medskip

\n {\bf Proof of Lemma 4.} Plainly it suffices to find $(\theta_i)^n_1$
positive and $U$ orthogonal so that $(\theta_i(AU)_{ij})$ is positive and
has constant, non-zero, diagonal.

Since the rows of $A$ are non-null, there is a constant $c>0$ so that if
$(\theta_i)^n_1$ is a sequence of positive numbers

$$\|(\theta_ia_{ij})\|_{C_1} \ge c\max_i \theta_i.$$

\n Since $\|(\theta_ia_{ij})\|_{C_1}$ is continuous as a function of
$(\theta_i)^n_1$, there is a sequence $(\theta_i)^n_1$ of positive numbers
which minimises

$$\|(\theta_ia_{ij})\|_{C_1}$$

\n subject to the condition $\prod\limits_i \theta_i = 1$. Let $H =
(h_{ij})$ be the positive square root of\hfil\break
$(\theta_i(AA^*)_{ij}\theta_j)$,
for this particular sequence, and note that there is an orthogonal matrix
$U$ for which

$$h_{ij} = \theta_i(AU)_{ij}, \quad 1\le i,j\le n.$$

Again, since $A$ has non-null rows, the diagonal entries of $H$ are
non-zero. For each $i$, let

$$\gamma_i = {1\over \sqrt{h_{ii}}} \left( \prod^n_{j=1}
\sqrt{h_{jj}}\right)^{1/n}.$$

Since $\prod\limits^n_i \gamma_i =1$, the matrix $(\gamma_i\theta_ia_{ij})$
has nuclear norm at least that of $(\theta_ia_{ij})$, the latter being
$\|H\|_{C_1}$. So

$$\eqalign{\|H\|_{C_1} &\le \|(\gamma_i\theta_ia_{ij})\|_{C_1}\cr
&=\|(\gamma_ih_{ij})\|_{C_1}\cr
&= \left( \prod^n_{k=1} \sqrt{h_{kk}}\right)^{1\over n} \left\|\left(
{1\over \sqrt{h_{ii}}} h_{ij}\right)\right\|_{C_1}\cr
&\le \sqrt{n} \|H\|^{1\over 2}_{C_1} \left(\prod_k
\sqrt{h_{kk}}\right)^{1\over n}}$$

\n by Lemma 6. So

$$\left( {1\over n} \sum_i h_{ii}\right)^{1\over 2} \le \left(\prod_k
h_{kk}\right)^{1\over 2n}$$

\n implying that the $h_{ii}$'s are all the same.$\hfill \square$\medskip

\n {\bf Proof of Theorem 2$'$.} The statement to be proved is that if $A =
(a_{ij})$ is a real $n\times n$ matrix with 1's on the diagonal and
$(m_i)^n_1$ is a sequence of reals, then there is a sequence
$(\lambda_j)^n_1$ of reals with

$$\sum_j \lambda^2_j \le {1\over n}$$

\n and, for every $i$

$$\left|\sum_j a_{ij}\lambda_j - m_i\right| \ge {1\over n}.$$

Using Lemma~4, choose a sequence $(\theta_j)^n_1$ of positive numbers and
an orthogonal matrix $U$, so that if

$$H = (\theta_i (AU)_{ij}),\eqno (2)$$

\n then $H$ is positive and has 1's on the diagonal.

By Lemma~3, there is a choice of signs $(\varepsilon_j)^n_1$ so that for
each $i$,

$$\left|\sum_j h_{ij}\varepsilon_j\theta_j - n \theta_i m_i\right| \ge
\theta_i.\eqno (3)$$

\n From (2) and (3), one has that for each $i$,

$$\left|\theta_i \sum_j (AU)_{ij}\varepsilon_j\theta_j - n \theta_i
m_i\right|\ge \theta_i,$$

\n and hence

$$\left|\sum_k a_{ik}\left({1\over n}\sum_j
u_{kj}\varepsilon_j\theta_j\right) - m_i\right| \ge {1\over n}.$$

\n For each $k$ set

$$\lambda_k = {1\over n} \sum_j u_{kj}\varepsilon_j \theta_j.$$

\n It remains to check that $\sum \lambda^2_k \le {1\over n}$. But
$\sum \lambda^2_k = {1\over n^2} \sum \theta ^2_j$ since $U$ is orthogonal
and so what is needed is

$$\sum_j \theta^2_j \le n.$$

From (2),

$$\theta_i a_{ij} = (HU^*)_{ij} \quad \hbox{for all}\quad i \ {\rm and} \
j$$

\n and so in particular, since $a_{ii} = 1$ for each $i$,

$$\theta_i = (HU^*) _{ii},\quad 1\le i \le n.$$

\n Now since $h_{ii} = 1$ for each $i$, Lemma~5 shows that

$$\eqalignno{\sum_i \theta^2_i &= \sum_i (HU^*)^2_{ii} = \sum_i
{(HU^*)^2_{ii}\over h_{ii}}\cr
&\le \sum_i h_{ii} = n.&\square}$$

\medskip

\n {\bf  \S 3. The infinite-dimensional case.}

Theorem~2 and weak$^*$-compactness immediately imply the following
``multiple Hahn-Banach'' theorem.\medskip

\n {\bf Theorem 7.} Let $(x_i)^\infty_1$ be a sequence of unit vectors in a
normed space $X, (m_i)^\infty_1$ a sequence of real numbers and
$(w_i)^\infty_1$ a sequence of non-negative reals with
$\sum\limits^\infty_1 w_i\le 1$. Then there is a linear functional $\phi$
of norm at most 1 in $X^*$ with

$$|\phi(x_i) - m_i| \ge w_i\quad \hbox{for every}\quad i.$$

For reflexive spaces, Theorem~1 follows immediately from Theorem~2. For
general spaces, Theorem~1 is a little more delicate. It can be regarded as
a quantitative strengthening of the Banach-Steinhaus theorem. If $\Phi$ is
an unbounded subset of the dual $X^*$ of a Banach space $X$, then there are
elements $\phi_1,\phi_2,\ldots$ of $\Phi$ with (say)

$$\sum^\infty_1 n\|\phi_n\|^{-1} < 1.$$

\n By Theorem~1, there is a point $x\in X$ of norm at most 1 so that for
each $n$,

$$\left|{\phi_n\over \|\phi_n\|}(x)\right| > n\|\phi_n\|^{-1},$$

\n i.e. $|\phi_n(x)| > n$.

To prove Theorem~1 it is necessary to examine the proof of Theorem~2 more
closely. The change of density argument in Section~1 and the proof in
Section~2 actually yield the following strong form of Theorem~2$'$.
Theorem~1 will be deduced from this.\medskip

\n {\bf Theorem 8.} Let $(a_{ij})$ be a real, $n\times n$ matrix with 1's
on the diagonal, $(m_i)^n_1$ a sequence of real numbers and $(w_i)^n_1$, a
sequence of positive numbers. Then there is a sequence $(\lambda_j)^n_1$
with

$$\sum_j w^{-1}_j \lambda^2_j \le \sum_j w_j$$

\n and, for every $i$,

$$\left|\sum_j a_{ij}\lambda_j - m_i\right| \ge w_i.\eqno \square$$

\n {\bf Proof of Theorem 1.} Suppose $(\phi_i)^\infty_1$ are unit
functionals
in $X^*, (m_i)^\infty_1$ are real numbers and $(w_i)^\infty_1$ are
non-negative numbers with $\sum\limits_i w_i<1$. The problem is to find a
point $x$ in the unit ball of $X$ with

$$|\phi_i(x) - m_i| > w_i \quad \hbox{for each}\quad i.$$

\n Choose a sequence $(v_i)^\infty_1$ with

$$v_i > w_i \ge 0 \quad \hbox{for each}\quad i$$

\n but

$$\sum_i v_i = 1 - \varepsilon < 1.$$

\n For each $i$, choose a point $x_i$ of norm at most 1 in $X$ with
$\phi_i(x_i) = 1-\varepsilon$. Applying Theorem~8 to the first $n$
functionals and vectors one obtains, for each $n$, a sequence
$(\lambda^{(n)}_j)^n_{j=1}$ satisfying

$$\sum^n_{j=i} v^{-1}_i(\lambda^{(n)}_j)^2 \le (1-\varepsilon)^{-2}
\sum^n_{j=i} v_j < (1-\varepsilon)^{-1}\eqno (5)$$

\n and for $1\le i \le n$

$$\left| \phi_i\left(\sum^n_{j=i} \lambda^{(n)}_jx_j\right) - m_i\right|
\ge v_i.$$

\n Regard $(\lambda^{(n)}_j)$ as an infinite sequence by filling out with
zeroes. From (5), for each $n$,

$$\eqalignno{\sum^\infty_{j=1} |\lambda^{(n)}_j| &\le \left(\sum^\infty
_{j=1} v_j\right)^{1/2} \left(\sum^\infty_{j=1} v^{-1}_j
(\lambda^{(n)}_j)^2 \right)^{1/2}&(6)\cr
&< (1-\varepsilon)^{1/2}(1-\varepsilon)^{-1/2}=1.}$$

\n Moreover, for each $m$ and $n$,

$$\sum^\infty_{j=m} |\lambda^{(n)}_j| \le \left(\sum^\infty_{j=m}
v_j\right)^{1/2} (1-\varepsilon)^{-1/2}.$$

\n Since the right-hand side $\to 0$ as $m\to \infty$, the sequences
$(\lambda^{(n)}_j)^\infty_{j=1}$ are uniformly summable, so the collection
has a (norm) limit point $(\lambda_j)^\infty_1$ (say) in $\ell_1$. From
(6), the point $x = \sum\limits _j \lambda _jx_j \in X$ has norm at most 1
and clearly

$$|\phi_i(x) - m_i| \ge v_i> w_i \quad \hbox{for every}\quad i.\eqno
\square$$

\vfill\eject

\n {\bf References}

\item{[A]} R. Alexander, A problem about lines and ovals, Amer. Math.
Monthly 75 (1968), 482-487.

\item{[B]} T. Bang, A solution of the ``Plank problem'', Proc. Amer. Math.
Soc. 2 (1951), 990-993.

\item{[D]} H. Davenport, A note on Diophantine approximation, Studies in
mathematical analysis and related topics (Stanford University Press, 1962),
77-81.

\item{[Ga]} R.J. Gardner, Relative width measures and the plank problem,
Pacific J. Math. 135 (1988), 299-312.

\item{[Gr]} H. Groemer, Coverings and packings by sequences of convex
sets.
In Discrete Geometry and Convexity, Annals of N.Y. Acad. Sci. Vol. 440
(1985).

\item{[R]} C.A. Rogers, Some problems in the geometry of convex sets, In
The Geometric Vein, Davis, Gr\"unbaum \& Scherk (eds.), Springer-Verlag,
New York (1981).

\end